\documentclass[11pt,reqno]{amsart}
\usepackage{amssymb,amsthm,amsmath}
\usepackage{floatrow}
\usepackage{mathabx}
\usepackage{float}
\usepackage{mathtools}
\usepackage{subfigure}
\usepackage{geometry}
\usepackage{lipsum}
\usepackage{caption} 
\usepackage{hyperref}
\usepackage{tikz}
\usepackage{array}
\usepackage{enumitem}
\usepackage{multirow}
\usepackage{enumitem}   
\usepackage{tikz-cd}

\usepackage{graphicx}

\divide
\oddsidemargin by 2 
\date{}
\makeatletter
\def\ps@pprintTitle{
  \let\@evenfoot\@oddfoot
}
\makeatother

\newtheorem{theorem}{Theorem}

\newtheorem{lemma}{Lemma}
\newtheorem{question}{Question}
\newtheorem{corollary}{Corollary}

\makeatletter
\newcommand{\namelabel}[1]{%
  \phantomsection
  \renewcommand{\@currentlabel}{#1}
  \label{#1}
}
\makeatother

\begin{document}
\title{Lifting closed curves to finite  covers of free groups
}

 \author{Deblina Das}
 \address{Department of Mathematics, Indian Institute of Technology Palakkad}
 \email{212114002@smail.iitpkd.ac.in, deblina099@gmail.com}
 \author{Arpan Kabiraj}
 \address{Department of Mathematics, Indian Institute of Technology Palakkad}
\email{arpaninto@iitpkd.ac.in}
\begin{abstract}
In this article, we show that given any integer $l\geq 2$, every closed curve $\gamma$ on the bouquet of $n$-circles $\Gamma$, admits a lift to a finite $l$-sheeted normal covering of $\Gamma$. Equivalently, identifying the free group $F_n$ of $n$ generators with the fundamental group of $\Gamma$, this statement asserts that $F_n$ is a union of $ l$-index normal subgroups for any $l\geq 2.$ The proof proceeds by explicitly constructing families of $l$-sheeted normal coverings of $\Gamma$, together with a characterization, in terms of necessary and sufficient conditions, of when a closed curve $\gamma$ on $\Gamma$ lifts to these covers.
\end{abstract}
\maketitle
\section{Introduction}

Let $\Gamma$ be the bouquet of $n$ circles, i.e. the wedge sum of $n$ circles $\bigvee_{i=1}^{n} S^1$ at the vertex $w$. The fundamental group $\pi_1(\Gamma)$ of $\Gamma$ is the free group $F_n$ with $n$ generators. We label these generators by $a_1,a_2,\ldots, a_n$, identifying each $a_i$ with the corresponding circle in $\Gamma$ at $w$.
A closed curve on $\Gamma$ based at $w$ is then naturally identified with a word in the free group $F_n$.
Let $H$ be a finite index subgroup of $F_n$ and $A$ be any generating set of $H$. In his seminal paper \cite{stallings1983topology}, Stallings gave an explicit construction of a finite sheeted cover of $\Gamma$ corresponding to the subgroup $H$, based on the generating set $A$. Our work is motivated by the following question, which may be viewed as dual to Stallings’ construction.
\begin{question}
Given any integer $l\geq 2$ and a closed curve  $\gamma$ on $\Gamma$ based at $w$, can we explicitly construct an $l$-sheeted normal covering of $\Gamma$ in which $\gamma$ admits a lift? In other words, given any integer $l\geq 2$ and $\gamma\in F_n$, can we explicitly construct an $l$-index normal subgroup of $F_n$ containing $\gamma$?
\end{question}

We answered the question affirmatively. 
To address Question $1$, we construct five families $M_u, M_{u,v}^k, M_{v,u}^k,N_{u,v}^k, N_{v,u}^k$ (see Section \ref{sec:pre} for details) of normal $l$-sheeted abelian coverings of $\Gamma$, where $l$ is an odd prime and two families of coverings for $l=2$. For each of these families of coverings, we provide necessary and sufficient conditions for a closed curve $\gamma$ on $\Gamma$ to admit a lift to those covers in terms of the word presentation of $\gamma$. We subsequently apply prime factorization together with an algebraic Lemma \ref{Lemma on integers}, to obtain a complete answer to Question $1$ for all integers $l\geq 2$. 
As a consequence, we establish the 
following theorem.
\begin{theorem}\label{main theorem}
    Let  $l \geq 2$ be any positive integer. Given any closed curve $\gamma$ on the bouquet of $n$ circles $\Gamma$ based at $w$, there exists an $l$-sheeted normal covering $\widetilde{\Gamma}$ of $\Gamma$ such that $\gamma$ lifts to $\widetilde{\Gamma}$. In other words, for any positive integer $l \geq 2$, $$\bigcup \big\{H\trianglelefteq F_n:[F_n:H]=l\big\}=F_n . $$
\end{theorem}

 Any compact, oriented surface $X$ with a non-empty boundary deformation retracts the bouquet of $n$ circles for some $n$. Using this 
association, as an application of Theorem \ref{main theorem}, we obtain the following result.
\begin{corollary}    
\label{Application on compact surface}
  For any positive integer $l \geq 2$ and any closed curve $\gamma$ on a compact oriented surface $X$ with a non-empty boundary, $\gamma$ lifts to some normal $l$-sheeted covering of $X$.  
\end{corollary}
  In view of the Corollary \ref{Application on compact surface}, it is a natural question to ask: does any closed curve $\gamma$ on a closed oriented surface $X$ lift to some normal $l$-sheeted covering of $X$ for a given positive integer $l (\geq 2)$? This conclusion indeed follows from our results. See Theorem \ref{closed surfaces}. For more related works in this context,  see
\cite{NeumannMR62122},\cite{MR2510838},\cite{bryce1997covering}.

The paper is organised as follows. In Section \ref{sec:pre}, we introduce the required notations and definitions and provide the construction of the families of finite sheeted coverings.  In section \ref{Index two subgroups of free groups}, we prove the Theorem \ref{main theorem} for $2$-sheeted covering spaces of $\Gamma$. Section \ref{section Lifting Criteria} establishes a necessary and sufficient condition for lifting a closed curve to these families of finite coverings. In section \ref{Index three subgroups of free groups}, we prove Theorem \ref{main theorem} for $p$-sheeted normal covering where $p$ is an odd prime, and later on we extend the conclusion for any positive integer $l\geq 2$. 

\section{Preliminaries}\label{sec:pre}
Let $\Gamma$ be the bouquet of $n$ circles, $\bigvee_{i=1}^{n} S^1$ with the vertex $w$. The fundamental group of $\Gamma$ is the free group $F_n$ with $n$ generators $a_1,a_2,\ldots,a_u,\ldots, a_v,\ldots  a_n$ corresponding to the circles on $\Gamma$ based at $w$. Given any closed curve $\gamma$ on $\Gamma$ based at $w$, we consider its homotopy class as an element of $\pi_1(\Gamma)$.  Abusing notation, we denote this element of $\pi_1(\Gamma)$  again by $\gamma$. 

With the above identification, a closed curve  $\gamma$ on $\Gamma$ based at $w$ is given by an unique word $x_1^{m_1}x_2^{m_2} \ldots x_r^{m_r}$ where $m_1,m_2\ldots m_r\in \mathbb{Z} $ and $$x_i \in \{a_1,a_2\ldots, a_n, a_1^{-1} ,a_2^{-1},\ldots a_n^{-1}\},\  x_i \neq x_{i+1}\ \text{and}\ x_i^{-1} \neq x_{i+1}.$$ Throughout the article, whenever we consider a closed curve $\gamma$ and its word representation, we use the above identification.

We associate different colors and directions of the generators so that $\Gamma$ becomes a directed colored graph with one vertex $w$ and $n$ directed colored edges $a_1,a_2,\ldots, a_n$.
 Recall any $l$-sheeted  covering space of $\Gamma$ is a connected directed colored graph with: i) $l$ vertices,
 ii) $n$ colors of the directed edges,
iii) the in-degree and the out-degree in any color equal to $1$, for any vertex and
    iv)  a chosen base vertex  (see \cite{FREE} and \cite{Hall_1949}). We follow this construction throughout the article to define finite-sheeted coverings of $\Gamma$.

Let $\gamma$ be a closed curve on $\Gamma$ based at $w$. Suppose  $u,v\in \{1,2,\ldots n\}$ such that $p_1,p_2,\ldots,p_{s_1}$ are powers of $a_u$; $q_1,q_2,\ldots,q_{s_2}$ are powers of $a_v$ in word representation of $\gamma$ where $p_i,q_j\in \{m_1,m_2,\ldots , m_r\}$ for all $i,j$.

We denote $o_{\gamma}(a_u)= \displaystyle\sum_{h=1}^{s_1} p_l$ and $o_{\gamma}(a_v)=\displaystyle\sum_{h=1}^{s_2} q_l$.\label{word representation} 


We denote the covering projection by $p$. In an $ l$-sheeted covering space, we fix an order of the vertices and denote them by $v_1,v_2,\ldots, v_l$, which are the preimages of vertex $w$ of $\Gamma$. Note that for each $j=1,2,\ldots,n$, the preimage $p^{-1}(a_j)$ has $l$ directed edges labeled as $a_j$-edges, joining the vertices $v_i$ where $i\in \{1,2,\ldots l\}$. 
\section{Index two subgroups of free groups}\label{Index two subgroups of free groups}
  We consider two families of $2$-sheeted coverings: $M_i$ and $M_{i,j}$, and show that any closed curve on $\Gamma$ lifts to one of these coverings. 
\begin{theorem}\label{Theorem free group 2}
    Any closed curve $\gamma$ on the bouquet of $n$ circles $\Gamma$ lifts to some $2$-sheeted  covering of $\Gamma$. In other words, the set $$\bigcup \big\{H\trianglelefteq F_n:[F_n:H]=2\big\}=F_n .$$
\end{theorem}
\begin{proof}

Let $\gamma$ be a closed curve on $\Gamma$ given by the word representation \ref{word representation}. We consider $a_i$ and $a_j$ for some $i,j\in \{1,2,\ldots,n\} $ and $i\neq j$.

  \textit{Type $M_i$:}\namelabel{$Type M_i$} 
         $p^{-1}(a_i)$ is connected such that we join $v_1\to v_2\to v_1$ by $a_i$-edges and the preimage of other generators has two connected components each as shown in Figure \ref{L_i,L_ij}(a). 
         In this case we observe if $o_{\gamma}(a_i)$ is even, $\gamma$ lifts to $M_i$, $i=1,2,\ldots,n$. 

   \textit{Type $M_{i,j}$:}\namelabel{$Type M_{i,j}^1$} 
          $p^{-1}(a_i)$ and $p^{-1}(a_j)$ are connected such that we join $v_1\to v_2\to v_1$ by the $a_i$-edges and $a_j$-edges, and the preimage of other generators has two connected components each as depicted in Figure \ref{L_i,L_ij}(b).  
         In this case we observe if $o_{\gamma}(a_i)+o_{\gamma}(a_j)$ is even, $\gamma $ lifts to $M_{i,j}$, $i,j\in \{1,2,\ldots,n\} $ and $i\neq j$. 

   \begin{figure}[h]
    \centering
   \subfigure[$2$-sheeted covering space $M_i$.]{\includegraphics[width=7.5 cm, height=4cm]{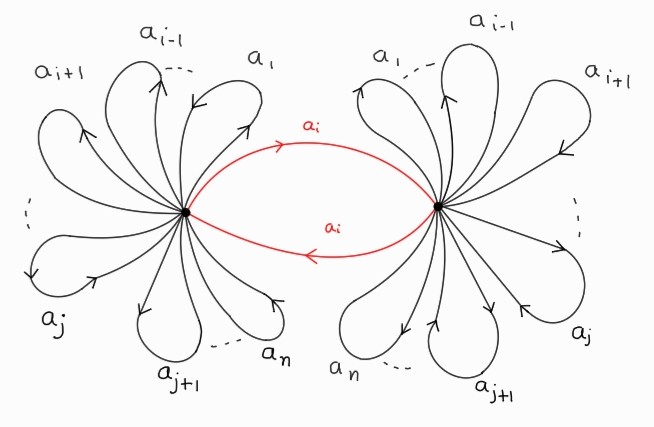}}  
    \subfigure[$2$-sheeted covering space $M_{i,j}$ , $i\neq j$.]{\includegraphics[width=7.5cm, height=4cm]{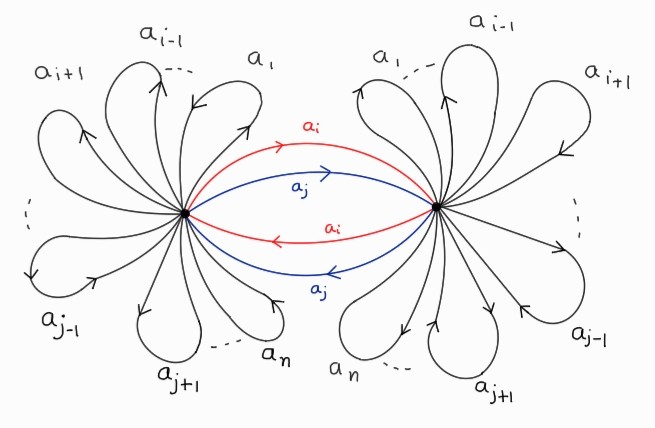}} 
    \caption{Two types of $2$-sheeted coverings of $\Gamma$.}
  \label{L_i,L_ij}

\end{figure}

    Note that it does not matter which one is the base point: this is consistent with the fact that any subgroup of index $2$ is normal.
    
    \begin{itemize}
        \item If $o_{\gamma}(a_i)+o_{\gamma}(a_j)$ is even, $\gamma$ lifts to $M_{i,j}$. Moreover, if $o_{\gamma}(a_i)$ and $o_{\gamma}(a_j)$ both even, $\gamma$ lifts to both $M_i$ and $M_j$ as well.
        \item If $o_{\gamma}(a_i)+o_{\gamma}(a_j)$ is odd, then one of $o_{\gamma}(a_i)$ and $o_{\gamma}(a_j)$ has to be even. Then $\gamma$ lifts to $M_i$ or $M_j$ depending on whether $o_{\gamma}(a_i) $ or $o_{\gamma}(a_j)$ is even, respectively.
    \end{itemize}
   
   Consequently, $\gamma$ lifts to some $2$-sheeted covering of $\Gamma$.
The second conclusion follows from the covering correspondence theorem \cite{Hatcher:478079}. 
\end{proof}

\section{Lifting Criteria}\label{section Lifting Criteria}
  In this section, we present explicit descriptions of the five families of coverings considered in our approach to the proof of Theorem \ref{main theorem}. 
  
    For any $u,v\in \{1,2,\ldots n\}$ with $u\neq v$ we focus on the following five families of abelian covers. 
Let $l$ is an odd prime, we write $l=2m-1$ where $m\geq 2$. In the following, $k$ denotes any positive integer between $1$
 and $m-1$, i.e. $1\leq k\leq m-1$. 

\begin{figure}[H]
    \centering
    \subfigure[$7$-sheeted covering space $M_{v,u}^2$ of $\Gamma$.]{  \includegraphics[width=7 cm, height=6.5cm]{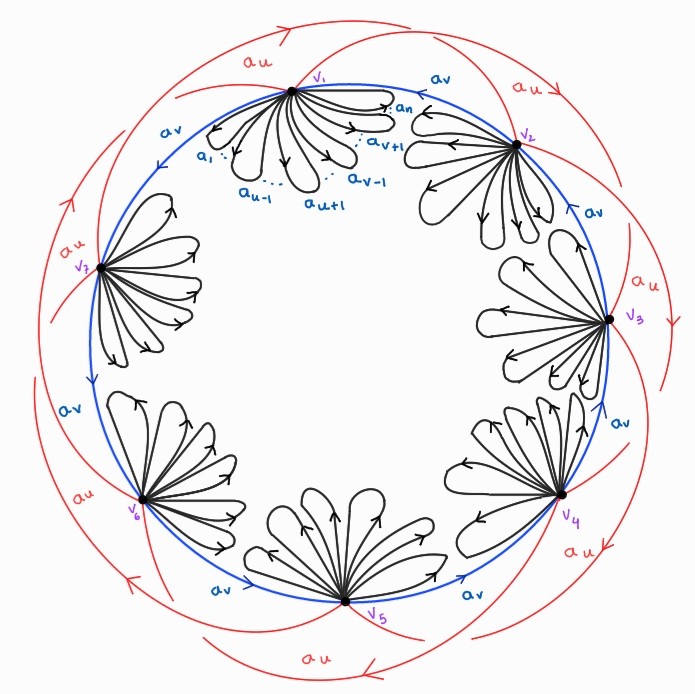}} 
    \subfigure[$7$-sheeted covering pace $N_{v,u}^2$ of $\Gamma$.]{    \includegraphics[width=7 cm, height=6.5cm]{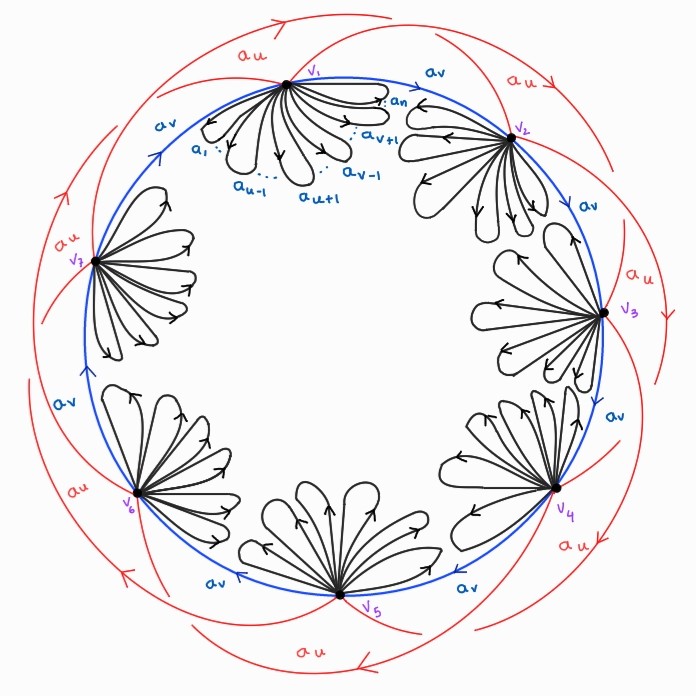} } 
     
    \caption{Five types of $7$-sheeted coverings of $\Gamma$.}
\label{Five types of 7-sheeted coverings of Gamma.}

\end{figure}

\textit{{\underline {Type $M_u$:}}} The inverse image $p^{-1}(a_u)$ is connected such that $v_i$ is connected to $v_{i\pmod{l}+1}$ by $a_u$-edges where $i=1,2,\ldots,l$ and $p^{-1}(a_j)$ has $l$ connected components each where $1\leq j\leq n, u\neq j$. 
       We refer to this family of covering as $M_u$ (see, Figure \ref{Five types of 7-sheeted coverings of Gamma.}(a)).
     
  \textit{{\underline { Type $M_{u,v}^k$:}}}  The inverse image $p^{-1}(a_u)$ is connected such that $v_i$ is connected to $v_{i\pmod{l}+1}$ by $a_u$-edges where $i=1,2,\ldots, l$. Also $p^{-1}(a_v)$ is connected such that $v_i$ and $v_{i+k-1\pmod{l}+1}$ being joined by $a_v$-edges. The directions of  $a_u$-edges and $a_v$-edges are different in the covering (if  $a_u$-edges are clockwise, then  $a_v$-edges are anticlockwise and vice-versa). The inverse image of other generators has $l$ connected components each. We refer to this family of covering as $M_{u,v}^k$ (see Figure \ref{Five types of 7-sheeted coverings of Gamma.}(b)).
  
  \textit{{\underline { Type $N_{u,v}^k$:}}}  The inverse image $p^{-1}(a_u)$ is connected such that $v_i$ is connected to $v_{i\pmod{l}+1}$ by  $a_u$-edges where $i=1,2,\ldots, l$. Also $p^{-1}(a_v)$ is connected such that $v_i$ and $v_{i+k-1\pmod{l}+1}$ being joined by  $a_v$-edges. The directions of  $ a_u$-edges and $a_v$-edges are the same in the covering (either both  $a_u$-edges and $a_v$-edges are clockwise or both  $a_u$-edges and $a_v$-edges are anticlockwise). The inverse image of other generators has $l$ connected components each. We refer to this family of covering as $N_{u,v}^k$ (see Figure \ref{Five types of 7-sheeted coverings of Gamma.}(c)).
 
  \textit{{\underline { Type $M_{v,u}^k$:}}}  The inverse image $p^{-1}(a_v)$ is connected such that $v_i$ is connected to $v_{i\pmod{l}+1}$ by  $a_v$-edges where $i=1,2,\ldots, l$. Also $p^{-1}(a_u)$ is connected such that $v_i$ and $v_{i+k-1\pmod{l}+1}$ being joined by  $a_u$-edges. The directions of  $a_u$-edges and $a_v$-edges are different in the covering (if  $a_u$-edges are clockwise, then  $ a_v$-edges are anticlockwise and vice-versa) . The inverse image of other generators has $l$ connected components each. We refer to this family of covering as  $M_{v,u}^k$ (see, Figure \ref{Five types of 7-sheeted coverings of Gamma.}(d)).

  \begin{figure}[h]
    \centering
   \subfigure[$7$-sheeted covering space $M_u$ of $\Gamma$.]{\includegraphics[width=7cm, height=6.5cm]{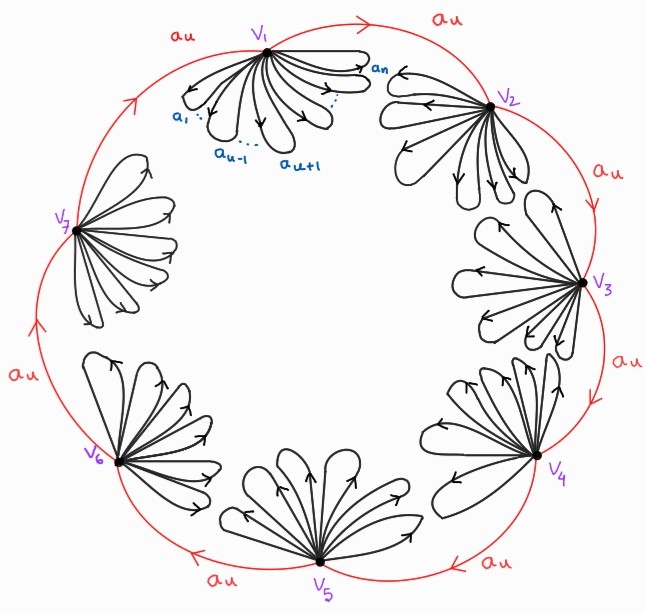} }  
   \end{figure}

    \begin{figure}[H]
    \centering
     \subfigure[$7$-sheeted covering space $M_{u,v}^3$ of $\Gamma$.]{  \includegraphics[width=7 cm, height=6.5cm]{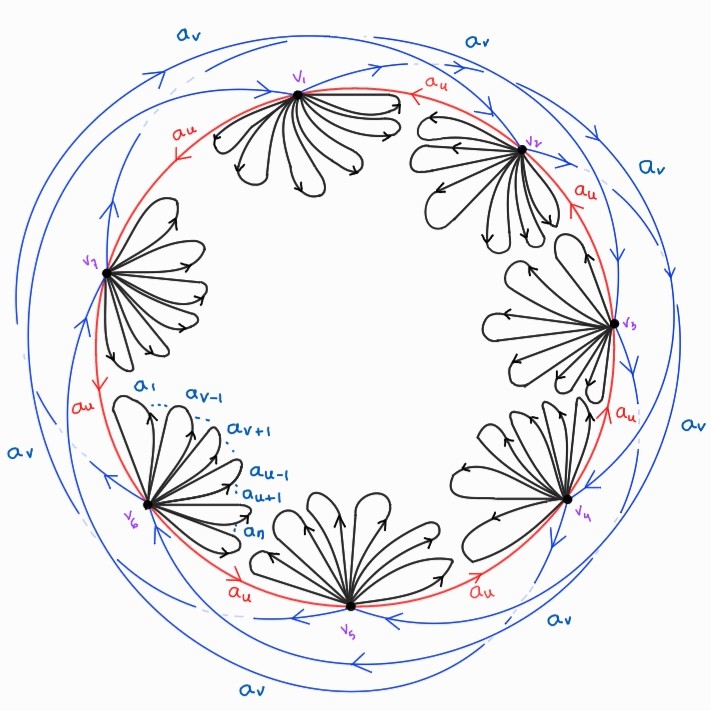}} 
   \subfigure[$7$-sheeted covering space $N_{u,v}^3$ of $\Gamma$.]{   \includegraphics[width=7 cm, height=6.5cm]{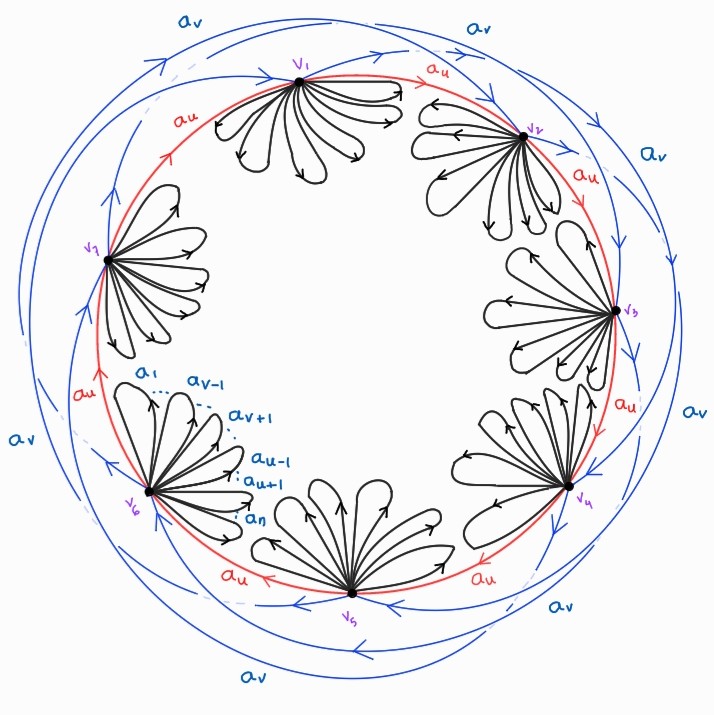} }  
   \end{figure} 
 \textit{{\underline { Type $N_{v,u}^k$:}}} The inverse image $p^{-1}(a_v)$ is connected such that $v_i$ is connected to $v_{i\pmod{l}+1}$ by  $a_v$-edges where $i=1,2,\ldots, l$. Also $p^{-1}(a_u)$ is connected such that $v_i$ and $v_{i+k-1\pmod{l}+1}$ being joined by  $a_u$-edges. The directions of  $a_u$-edges and $a_v$-edges are the same in the covering (either both  $a_u$-edges and $a_v$-edges are clockwise or both  $a_u$-edges and $a_v$-edges are anticlockwise). The inverse image of other generators has $l$ connected components each. We refer to this family of covering as $N_{v,u}^k$ (see Figure \ref{Five types of 7-sheeted coverings of Gamma.}(e)). 
 
   Now we are ready to obtain necessary and sufficient conditions for a curve to lift to the families of normal $l$-sheeted coverings of types $M_u$, $M_{u,v}^k$, $N_{u,v}^k$, $M_{v,u}^k$ and $N_{v,u}^k$ of $\Gamma$ where $l$ is odd prime. Note that, 
   it does not matter which one is the base point: this is consistent with the fact that all these covers are normal.

 \begin{lemma}\label{necessary sufficient condition lemma} Let $l$ be any odd prime and $l=2m-1$, where $m\in \mathbb{N}_{\geq 2 }$. For each $k$ where $1\leq k\leq m-1$, and for any $u,v\in \{1,2,\ldots,n\}$, with $u\neq v$,  
         \begin{enumerate}
     \item a closed curve $\gamma$ on $\Gamma$ lifts to $M_u$ if and only if $o_{\gamma}(a_u)\equiv0 \pmod{l}$.
     \item a closed curve $\gamma$ on $\Gamma$ lifts  to $M_{u,v}^k$ if and only if $o_{\gamma}(a_u)\equiv k o_{\gamma}(a_v) \pmod{l}$.
      \item a closed curve $\gamma$ on $\Gamma$ lifts to $N_{u,v}^k$ if and only if $k o_{\gamma}(a_v)+ o_{\gamma}(a_u)\equiv0 \pmod{l}$.
     \item a closed curve $\gamma$ on $\Gamma$ lifts  to $M_{v,u}^k$ if and only if $k o_{\gamma}(a_u)\equiv o_{\gamma}(a_v) \pmod{l}$.
    
     \item a closed curve $\gamma$ on $\Gamma$ lifts to $N_{v,u}^k$ if and only if $o_{\gamma}(a_v)+k o_{\gamma}(a_u)\equiv0 \pmod{l}$.

 \end{enumerate}
   \end{lemma}
   \begin{proof}
      \textit{(1)} Let $u\in \{1,2,\ldots,n\}$. From the Figure \ref{Five types of 7-sheeted coverings of Gamma.}(a) that the liftability condition of the closed curve  $\gamma$ on $M_u$ depends only on the occurrence of the generator $a_u$ in the word representation of $\gamma$; the other generators contribute no condition to that. The first claim immediately follows.
      
       Let $u,v\in \{1,2,\ldots,n\}$ such that $u\neq v$ and fix a $k$ such that $1\leq k\leq {m-1}$. From the Figure \ref{Five types of 7-sheeted coverings of Gamma.}, observe that the liftability condition of a closed curve  $\gamma$ to $M_{u,v}^k$  $N_{u,v}^k$, $M_{v,u}^k$ and $N_{v,u}^k$ depends only on the occurrence of the generators $a_u$ and $a_v$ in the word representation of $\gamma$ and independent of other generators, i.e., the other generators contribute no condition to that. Therefore, to check the necessary and sufficient conditions of liftability of a curve $\gamma$ to these covers, it is enough to consider curves $\gamma$ whose word representation consists of only the powers of $a_u$ and $a_v$. Without loss of generality, we consider a curve whose word representation starts with some power of $a_v$. \\
      \textit{(2)} On the cover $M_{u,v}^k$ , if we traverse along one $a_v$-edge, we have to traverse back along $k$ many  $a_u$-edges to form a closed curve. Thus, the smallest subword in $\gamma$ that forms a closed curve which
       lifts to $M_{u,v}^k$  is of the form  
 $a_v^{h\pmod{l}}a_u^{hk\pmod{l}}$ 
    where $h=0,1,2,\ldots, l-1$.
         The claim now follows from the induction on word length.
      Similar argument works for case \textit{(4)}.\\
      \textit{(3)}  On the cover $N_{u,v}^k$, if we traverse along one $a_v$-edge, we have to traverse back along $k$ many $a_u^{-1}$-edges to form a closed curve. 
      Thus, the smallest subword in $\gamma$ that forms a closed curve which
       lifts to $N_{u,v}^k$  is of the form 
       $a_v^{h\pmod{l}}a_u^{-hk\pmod{l}}$ 
    where $h=0,1,2,\ldots, l-1$.  The claim now follows from the induction on word length.
Similar argument works for case \textit{(5)}.

This completes the proof.
   \end{proof}
    \section{Finite index subgroups of free groups}\label{Index three subgroups of free groups}
In this section, we establish an algebraic lemma and, together with Lemma \ref{necessary sufficient condition lemma}, use it to complete the proof of Theorem \ref{main theorem}.
\begin{lemma}\label{Lemma on integers}
    Let $p$ be an odd prime and $p=2m-1$. For any pair of integers $a,b$, one of the following conditions is satisfied.
    \begin{enumerate}
        \item $ a\equiv 0 \pmod{p}.$
        \item $b\equiv 0 \pmod{p}. $
        \item $ax+b\equiv 0\pmod{p}$ for some $x$ such that $1\leq x\leq m-1.$
        \item $ax\equiv b \pmod{p}$ for some $x$ such that $1\leq x\leq m-1. $
    \end{enumerate}
    
\end{lemma}
\begin{proof}
    Consider any pair of integers $a,b\in \mathbb{Z}.$ If $a\equiv 0 \pmod{p}$, then $\textit{(1)}$ is satisfied.
    If $b\equiv 0 \pmod{p}$ then $\textit{(2)}$ is satisfied. If $a \not \equiv 0 \pmod{p}$ and $b\not \equiv 0 \pmod{p}$, the $a^{-1}$ and $b^{-1}$ both exist in the integer modulo $p$ field $\mathbb{Z}_p$. Let $k$ be a primitive solution of $x \equiv -a^{-1}b\pmod {p}$. Note that $k$ is unique as $gcd(a,p)=1$, . If $k\leq m-1$, then the condition $\textit{(3)}$ holds. If $k>m-1$ we claim that $p-k$ satisfies $ax\equiv b \pmod{p}.$ As, $$a(p-k)\equiv -ak \pmod {p} \equiv -a (-a^{-1}b) \pmod{p}\equiv b \pmod{p}.  $$
    Since $k>m-1$, we have $0<p-k<m$ , i.e., $p-k\leq m-1$. Therefore, $\textit{(4)}$ is satisfied. This completes the proof.
\end{proof}
\noindent\textbf{Proof of Theorem \ref{main theorem}:}\begin{proof}
 We consider two cases:  Case 1: when the index $l$ is  prime and Case 2: when the index $l$ is composite.
 
\textit{Case 1:} If $l=2$, the conclusion follows from Theorem \ref{Theorem free group 2} in section \ref{Index two subgroups of free groups}. 

    Let $ l=p$ be an odd prime number and $\gamma$ be a closed curve on $\Gamma$ based at $w$ given by the word representation described in \ref{word representation}.
Let, $u,v\in \{1,2,\ldots n\}$ where $u\neq v$.  
Now $o_{\gamma}(a_u), o_{\gamma}(a_v) \in \mathbb{Z}$. 
Using Lemma \ref{Lemma on integers}, one of the following conditions holds.
\begin{enumerate}
    \item $o_{\gamma}(a_u)\equiv 0 \pmod{p}$. In this case $\gamma$ lifts to $M_u$ using Lemma \ref{necessary sufficient condition lemma}.
    \item $o_{\gamma}(a_v)\equiv 0 \pmod{p}$. If this case $\gamma$ lifts to $M_v$ using Lemma \ref{necessary sufficient condition lemma}.
    \item Either of one holds, \begin{enumerate}[label=(\roman*)] 
        \item $o_{\gamma}(a_u)+ko_{\gamma}(a_v)\equiv 0 \pmod {p} $ for some $k$ where $1\leq k\leq m-1$. In this case $\gamma$ lifts to $N_{u,v}^k$ for some $k$ where $1\leq k\leq m-1$ using Lemma \ref{necessary sufficient condition lemma}.
        \item $ko_{\gamma}(a_u)+o_{\gamma}(a_v)\equiv 0 \pmod{p} $ for some $k$ where $1\leq k\leq m-1$. In this happens, $\gamma$ lifts to $N_{v,u}^k$ for some $k$ where $1\leq k\leq m-1$ using Lemma \ref{necessary sufficient condition lemma}.
    \end{enumerate} 
    \item Either of one holds, \begin{enumerate}[label=(\roman*)]
    \item $ko_{\gamma}(a_v)\equiv o_{\gamma}(a_u) \pmod{p}$ for some $k$ where $1\leq k\leq m-1$. If this happens, $\gamma$ lifts to $M_{u,v}^k$ for some $k$ where $1\leq k\leq m-1$ using Lemma \ref{necessary sufficient condition lemma}.
        \item $o_{\gamma}(a_v)\equiv ko_{\gamma}(a_u) \pmod{p}$ for some $k$ where $1\leq k\leq m-1$. If this happens, $\gamma$ lifts to $M_{v,u}^k$ for some $k$ where $1\leq k\leq m-1$ using Lemma \ref{necessary sufficient condition lemma}.
    \end{enumerate} 
\end{enumerate} 
Therefore, $\gamma$ lifts to some $p$-sheeted normal covering of $\Gamma$ where $p$ is an odd prime. 

\textit{Case 2:} Let $l$ be composite and $l=p_1^{r_1}p_2^{r_2}\ldots p_k^{r_k}$ be the prime decomposition of $l$ where $p_1,p_2,\ldots p_k$ are distinct prime numbers and $r_i (1\leq i\leq k)$ are positive integers. Recall that if we have a tower of subgroups, $H\leq K\leq G$, then their indices satisfy the relation  $[G:H]= [G:K][K:H]$. Combining this with the fact that every subgroup of a free group is free, successive application of \textit{Case 1} concludes the proof in this case.

The second conclusion follows from the covering correspondence theorem \cite{Hatcher:478079}.
\end{proof}
 \begin{theorem}\label{closed surfaces}
  Any closed curve $\gamma$ on a closed oriented surface $X$ lifts to some normal finite sheeted covering of $X$.
 
  \end{theorem}
  \begin{proof}
      Let $S_g$ be a compact surface of genus $g$  with generators $c_1,d_1,\ldots, c_g, d_g$ and $Q$ be a regular neighbourhood of the generators. Let $\gamma \in \pi_1(S_g)$ be a non-trivial element. There is a canonical projection $i_*: F_{2g} \to \pi_1(S_g)$ as $\pi_1(S_g)=F_{2g}/N$ where $N=\langle\langle R \rangle\rangle$ and $R=[c_1,d_1][c_2,d_2]\ldots[a_g,d_g]$.
    Since $\gamma \in \pi_1(S_g),$ by Corollary \ref{Application on compact surface} there exists a finite sheeted normal covering $p:\Tilde{Q}\to Q$ such that preimage of $\gamma$ (by $i_*$) lifts, i.e., lies in $H= p_*(\pi_1(\Tilde{Q}))\subset \pi_1(Q)=F_{2g}$. 
      Since $H=p_*(\pi_1(\Tilde{Q}))$ is a normal subgroup of $F_{2g}$ and $i_*$ is surjective, $i_*(H)$ is normal subgroup of $\pi_1(S_g)$. We claim $i_*(H)$ is finite index subgroup of $\pi_1(S_g)$. We  define a map $\chi : \pi_1{(Q)}/ H\to \pi_1(S_g)/i_*(H)$ by $\Psi(xH)=i_*(x)i_*(H) $ which is surjective. Since $\pi_1{(Q)}/ H$ is finite, $\pi_1(S_g)/i_*(H)$ is finite. Using covering correspondence theorem, there exists a finite sheeted normal covering $q: \Tilde{S_g}\to S_g$ such that $q_*(\pi_1(\Tilde{S_g}))=i_*(H)$. Hence  $\gamma$  lies in $ q_*(\pi_1(\Tilde{S_g}))$, i.e., lifts to $\Tilde{S_g}$ by $q$.
      \end{proof}
\begin{question}
    Given a closed curve $\gamma$ on $\Gamma$ based at $w$, can one determine criteria characterizing all finite sheeted coverings of $\Gamma$ to which $\gamma$ lifts, expressed in terms of the word presentation of $\gamma$?
\end{question}

\bibliographystyle{plain}
\bibliography{bibfile}

\end{document}